\documentstyle[12pt]{article}
\newfont{\frak}{eufm10 scaled\magstep1}
\newfont{\Bbb}{msbm10 scaled\magstep1}
\let\al\alpha
\let\be\beta
\let\ga\gamma
\let\de\delta

\let\ka\kappa

\let\De\Delta
\let\ep\varepsilon
\newcommand{\ZZ}{\hbox{\Bbb Z}}
\newcommand{\CC}{\hbox{\Bbb C}}

\def\taa{{\tilde A}_1}
\def\ta{{\tilde A}_0}
\def\tb{{\tilde B}_0}
\def\tc{{\tilde C}_0}

\newcommand{\FSA}{{\cal A}}

\newcommand{\FSD}{{\cal D}}

\newcommand{\FSU}{{\cal U}}

\let\iy\infty

\newcommand{\tfrac}[2]{{\textstyle\frac{#1}{#2}}}

\newcommand{\Proof}{\noindent{\bf Proof}\ \ }

\newtheorem{theorem}{Theorem}[section]

\newtheorem{lemma}[theorem]{Lemma}
\newtheorem{Remark}[theorem]{Remark}

\begin{document}
\vskip 50mm
\title{${}_3F{}_2(1)$ hypergeometric function and \\
quadratic $R$-matrix algebra }
\vskip 45mm
\author{Vadim B. Kuznetsov${}\sp {\,1,2}$}
\footnotetext[1]{The author was supported by the
Netherlands Organisation for Scientific Research
(NWO) under the Project \# 611--306--540.}
\footnotetext[2]{On leave from Department of Mathematical and Computational
Physics, Institute of Physics, St.Petersburg University, St.Petersburg
198904, Russia.}
\vskip 10mm
\maketitle
{\it Department of Mathematics and Computer Science,
University of Amsterdam,
\par Plantage Muidergracht 24,
1018 TV Amsterdam, The Netherlands
\par e-mail: {\tt vadim@fwi.uva.nl}}
\vskip 20mm
\begin{abstract}\noindent
We construct a class of representations
of the quadratic $R$-matrix algebra given by the reflection equation
with the spectral parameter,
$$
R{\,}(u-v)\,T^{(1)}(u)\,R{\,}(u+v)\,T^{(2)}(v)=
T^{(2)}(v)\,R{\,}(u+v)\,T^{(1)}(u)\,R{\,}(u-v),
$$
in terms of certain ordinary difference operators.
These operators turn out to act as parameter shifting operators on
the ${}_3F{}_2(1)$ hypergeometric function and its limit cases and on classical
orthogonal polynomials. The relationship with the factorisation method
will be discussed.
\end{abstract}
\vskip 20mm
\noindent
Mathematical preprint series, University of Amsterdam, Report 94-21,
October 11, 1994

\noindent hep-th/yymmnnn
\vskip 10mm
{\bf Key words:} quadratic $R$-matrix algebra, reflection equation,
hypergeometric \par\noindent functions, classical orthogonal polynomials,
Hahn polynomials,
recurrence \par\noindent relations, contiguous function relations
\vskip 10mm\noindent
\vskip 10mm
\pagebreak
%
%
\section{Introduction}
\setcounter{equation}{0}
\noindent
Let $V$ be a complex vector space. The quantum $R$-matrix
is a meromorphic operator-valued function $R:\CC\rightarrow$End$(V\otimes V)$
which satisfies the quantum Yang-Baxter equation of the form
\cite{ba82,fa94,sk91}
$$
R^{(12)}(u-v)R^{(13)}(u-w)R^{(23)}(v-w)=
R^{(23)}(v-w)R^{(13)}(u-w)R^{(12)}(u-v).
$$
Let us fix the following solution of this equation
\begin{equation}
R(u)=u+\ka P, \qquad u,\ka\in\CC,
\qquad P(x\otimes y)=y\otimes x,\quad x,y\in V.
\label{nn}\end{equation}
Consider $V=\CC^2$ then in the standard basis
we have the following $4 \times 4$ $R$-matrix
\begin{equation}
R(u)=\left(\matrix{
u+\ka&0&0&0\cr 0&u&\ka&0\cr 0&\ka&u&0\cr 0&0&0&u+\ka}\right).
\label{o1}\end{equation}
Introduce a $2\times 2$ matrix
\begin{equation}
U(u)=\left(\matrix{A(u)&B(u)\cr C(u)&D(u)}\right)
\label{o2}\end{equation}
with a priori non-commuting entries depending on a so-called
{\em spectral parameter} $u$ which is arbitrary complex.
The {\em QISM II algebra} \cite{kk93},
or the algebra given by the {\em reflection equation},
is the algebra generated by the
matrix elements of $U(u)$ for all $u$ subject to a quadratic
relation involving two $R$-matrices \cite{ks92,sk88}:
\begin{equation}
R(u-v)U^{(1)}(u)R(u+v-\ka)U^{(2)}(v)=
U^{(2)}(v)R(u+v-\ka)U^{(1)}(u)R(u-v).
\label{o4}\end{equation}
Here we use the notation
$U^{(1)}(u)=U(u)\otimes I$, \,$U^{(2)}(v)=I\otimes U(v)$.
The QISM II algebra, or the {\em $U$-algebra}
is an example of a {\em quadratic $R$-matrix
algebra}. {}From now on we assume that $\ka=1$. For $\ka\ne0$
in (\ref{o1}) this means no loss of generality.

In the present article we construct
a class of representations of very simple type
($U$-operators of rank 1) of the $U$-algebra.
We require moreover a certain symmetry property (unitarity)
for an $U$-operator.
We will consider $U$-operators (\ref{o2}) for which certain matrix elements
will be realized in terms of ordinary difference operators acting
as parameter shifting operators on the ${}_3F{}_2(1)$
generalized hypergeometric function
and its limit cases. Specialization then yields shift operator actions on
classical orthogonal polynomials (more specifically, on Hahn polynomials
which are generic ones for the representations under consideration).
We will also point out a close connection of our results
with the factorization method for second order
difference equations.

Our operators will act on special functions $F(u)$
which appear for each $u\in\CC$ as a solution of the equation
\begin{equation}
C(u)\,F(u)=0,
\nonumber\end{equation}
i.e., as functions annihilated by one of the two off-diagonal
elements (always chosen to be $C(u)$)
of an $U$-operator. The operators $A(u)$ and $-A(-u)$
then give the shifting of the parameter $u$ by $\pm 1$,
respectively:
\begin{equation}
A(u)F(u)=\Delta_-(u-\tfrac12)F(u-1),\quad
-A(-u)F(u)=\Delta_+(u+\tfrac12)F(u+1).
\nonumber\end{equation}
Here the $\Delta_\pm(u)$ are certain scalars
depending on $u$ which factorize the quantum determinant $\Delta(u)$
of an $U$-operator:
\begin{equation}
\Delta(u)=\Delta_+(u)\Delta_-(u).
\nonumber\end{equation}
The quantum determinant for the $U$-algebra is a certain quadratic
expression in the
generators with the property that it is the generating function for the
center of the algebra.
So, in an irreducible representation it is,
under suitable assumptions, scalar for each $u$.

We would like to thank Tom H. Koornwinder for many helpful discussions
during preparation of this paper.
%
%
\section{More about the QISM II algebra}
\setcounter{equation}{0}
\noindent
We will always assume the following relations
(symmetry property when changing the sign of $u$) for our
$U$-operators:
\begin{eqnarray}
-A(-u)&=&D(u)-(A(u)+D(u))/(2u+1),\nonumber\\
-D(-u)&=&A(u)-(A(u)+D(u))/(2u+1),
\nonumber\\
B(-u)&=&B(u),\quad C(-u)=C(u).
\label{141}\end{eqnarray}
Note that the second equality is implied by the first.
The equations (\ref{141}) can be rephrased as the unitarity property
$U^{-1}(-u+1)\sim U(u)$ (\cite{sk88}).

The {\em quantum determinant} for the $U$-algebra is defined as follows.
\begin{eqnarray}
\Delta(u)
&=&-D(-u+\tfrac12)D(u+\tfrac12)-C(u-\tfrac12)B(u+\tfrac12)
\nonumber\\
&=&-A(-u+\tfrac12)A(u+\tfrac12)-B(u-\tfrac12)C(u+\tfrac12)
\nonumber\\
&=&-D(u+\tfrac12)D(-u+\tfrac12)-C(u+\tfrac12)B(u-\tfrac12)
\nonumber\\
&=&-A(u+\tfrac12)A(-u+\tfrac12)-B(u+\tfrac12)C(u-\tfrac12).
\label{o10}\end{eqnarray}
The quantum determinant is the generating
function for the center of the algebra \cite{sk88}.

Relation (\ref{o4}) can be rewritten in the following
extended form as commutators between the algebra generators
$A(u),B(u),C(u)$, and $D(u)$.
\begin{eqnarray}
{[B,B]}&=&[C,C]=0,\label{214}\\
{[A,A]}&=&- (BC-{\tilde {BC}})/(u+v),\label{215}\\
{[D,D]}&=&- (CB-{\tilde {CB}})/(u+v),\label{216}\\
{[A,B]}&=&- (AB-{\tilde {AB}})/(u-v)-
(AB+BD)/(u+v-1)\nonumber\\
&&-(AB+BD-{\tilde {AB}}-{\tilde {BD}})/
(u-v)/(u+v-1),\label{217}\\
{[B,A]}&=&- (BA-{\tilde {BA}})/(u-v)+
({\tilde {AB}}+{\tilde {BD}})/(u+v-1),\label{218}\\
{[A,C]}&=&- (CA-{\tilde {CA}})/(u-v)+
({\tilde {CA}}+{\tilde {DC}})/(u+v-1)\nonumber\\
&&-(CA+DC-{\tilde {CA}}-{\tilde {DC}})/
(u-v)/(u+v-1),\label{219}\\
{[C,A]}&=&- (AC-{\tilde {AC}})/(u-v)-
(CA+DC)/(u+v-1),\label{220}\\
{[D,B]}&=&- (BD-{\tilde {BD}})/(u-v)+
({\tilde {AB}}+{\tilde {BD}})/(u+v-1)\nonumber\\
&&-(AB+BD-{\tilde {AB}}-{\tilde {BD}})/
(u-v)/(u+v-1),\label{221}\\
{[B,D]}&=&- (DB-{\tilde {DB}})/(u-v)-
(AB+BD)/(u+v-1),\label{222}\\
{[D,C]}&=&- (DC-{\tilde {DC}})/(u-v)-
(CA+DC)/(u+v-1)\nonumber\\
&&-(CA+DC-{\tilde {CA}}-{\tilde {DC}})/
(u-v)/(u+v-1),\label{223}\\
{[C,D]}&=&- (CD-{\tilde {CD}})/(u-v)+
({\tilde {CA}}+{\tilde {DC}})/(u+v-1),\label{224}\\
{[A,D]}&=&- (CB-{\tilde {CB}})(u+v+1)/(u^2-v^2),\label{225}\\
{[D,A]}&=&- (BC-{\tilde {BC}})(u+v+1)/(u^2-v^2),\label{226}\\
{[B,C]}&=&- (DA-{\tilde {DA}})(u+v-1)/(u^2-v^2)\nonumber\\
&&- (AA-{\tilde {DD}})/(u+v),\label{227}\\
{[C,B]}&=&- (AD-{\tilde {AD}})(u+v-1)/(u^2-v^2)\nonumber\\
&&- (DD-{\tilde {AA}})/(u+v).\label{228}
\end{eqnarray}
Here we use for brevity the notations:
$[A,B]$ means the commutator $[A(u),B(v)]$, where the first parameter
is $u$ and the second one is $v$; $DA$ stands for the noncommutative
operator product $D(u)A(v)$; and ${\tilde {DA}}$ signifies $D(v)A(u)$
(where $v$ is the first parameter), and so on.

The following Theorem was proved in \cite{kk93}.
\begin{theorem}\label{t7}
Let $W$ be a complex vector space on which the $U$-algebra
acts by an algebra representation. Suppose ${\cal D}$ is a subset
of $\CC$ of the form $\{u_0+m \mid m\in\ZZ,\allowbreak j_-<m<j_+\}$,
where $u_0\in\CC$ and $j_\pm=\pm\iy$ or integer, such that

(i) $\{w\in W\mid C(u)w=0\}$ is 1-dimensional for any $u\in {\cal D}$,

(ii) if $u\in {\cal D}$, $0\neq w\in W$ and $C(u)w=0$ then
\begin{equation}
A(u)w\quad \left\{\quad \matrix{
\neq 0,& u\neq u_0+j_-+1,\cr
 =0, & u= u_0+j_-+1,}\right.
\nonumber\end{equation}
\begin{equation}
-A(-u)w\quad \left\{\quad \matrix{
\neq 0, & u\neq u_0+j_+-1,\cr
 =0,& u= u_0+j_+-1.}\right.
\nonumber\end{equation}
For each $u\in {\cal D}$ choose $0\neq F(u)\in W$ such that $C(u)F(u)=0$.
Then
\begin{eqnarray}
{A(u)F(u)}&=&\Delta_-(u-\tfrac12)F(u-1),\qquad
u\in {\cal D},\quad u\neq u_0+j_-+1,\nonumber\\
{-A(-u)F(u)}&=&\Delta_+(u+\tfrac12)F(u+1),\qquad
u\in {\cal D},\quad u\neq u_0+j_+-1,
\nonumber\end{eqnarray}
for certain scalar functions $\Delta_\pm (u\pm\tfrac12)$.
Furthermore, the operator
$\Delta(u)$, when acting on ${\rm Span}\,\{F(v)\mid  v\in {\cal D}\}$,
is scalar for $u\in\FSD\pm\tfrac12$ and it satisfies
\begin{equation}
\Delta (u)=\quad \left\{\quad \matrix{
\Delta_+(u)\Delta_-(u)\,,& u\in\FSD\pm\tfrac12\,,&
u\neq u_0+j_\pm\mp\tfrac12\,,\cr
0\,,&\qquad u=u_0+j_\pm\mp\tfrac12\,.&}\right.
\nonumber\end{equation}
\end{theorem}
\vskip 2mm
%
%
\section{The function ${}_3F{}_2(1)$}
\setcounter{equation}{0}
\noindent
Let us consider the following generalized hypergeometric function:
$$
F\equiv{}_3F{}_2(1)=
{}_3F_2(a,b,c\,;d,e;1)=\sum_{n=0}^\iy{(a)_n(b)_n(c)_n\over(d)_n(e)_nn!}\,,
$$
$$
{\rm Re}\,(a+b+c-d-e-1)<-1,\qquad d,e\ne0,-1,-2,\ldots,
$$
$$
(\al)_n=\al(\al+1)\dots(\al+n-1),\qquad (\al)_0=1.
$$
If in this function one, and only one, of the variables is
increased or decreased by unity, the resultant function is said
to be contiguous to the $F$ above. We will use the
notation $\De_\al^\pm$ for the operators defined by the following:
$$
{\De_a^+F}=F(a+)-F,\qquad \De_d^-F=F(d-)-F.
$$
There are 45 relations each expressing $F$ linearly in terms of
two of its 10 contiguous functions (see, for instance, \cite{ra45}).
In what follows we will write down a full list of those 45 relations
using the operators $\De_\al^\pm$. We are giving operator's
equalities but it is implied that they are true only when acting on the $F$.
%
%
\subsection{Contiguous function relations}
The first 10 simple relations have the form
\begin{eqnarray}
&&\quad \al\De_\al^+=\be\De_\be^+,\qquad\al,\be\in
\{a,b,c\},\label{1}\\
&&\quad \al\De_\al^+=(\de-1)\De_\de^-,\qquad
\al\in\{a,b,c\},\quad \de\in\{d,e\},\label{2}\\
&&\quad(d-1)\De_d^-=(e-1)\De_e^-.\label{3}\end{eqnarray}
The next 10 relations have the following form. The triple
$(\al,\be,\ga)$ in the formulas below means a cyclic
permutation of $(a,b,c)$ while pair $(\de,\ep)$ means
that of $(d,e)$.
\begin{eqnarray}
&&(\al-d)(\al-e)\De_\al^--\al\ga=
(\be-d)(\be-e)\De_\be^--\be\ga,\label{15}\\
&&\qquad\qquad(\al-\ep)\De_\al^-=\tfrac{(\be-\de)(\ga-\de)}{\de}\De_\de^++
\tfrac{\be\ga}{\de},\label{12}\\
&&\tfrac{(a-d)(b-d)(c-d)}{d}\De_d^++\tfrac{abc}{d}=
\tfrac{(a-e)(b-e)(c-e)}{e}\De_e^++\tfrac{abc}{e}.\label{11}\end{eqnarray}

One can supply 25 more equalities using the above 20.
In the next 6 formulas the triple $(\al,\be,\ga)$ means
any permutation of $(a,b,c)$.
\begin{eqnarray}
&&(\al-d)(\al-e)\De_\al^-+\be(a+b+c+1-d-e)\De_\be^++
\be\ga=0.\label{z21}\end{eqnarray}
In the next 6 formulas the triple $(\al,\be,\ga)$ means
a cyclic permutation of $(a,b,c)$ and $\de\in\{d,e\}$.
\begin{eqnarray}
&&(\al-d)(\al-e)\De_\al^-+(\de-1)(a+b+c+1-d-e)\De_\de^-+
\be\ga=0.\label{z41}\end{eqnarray}
In the next 6 formulas $\al\in\{a,b,c\}$ and $\de\in\{d,e\}$.
\begin{eqnarray}
&&\tfrac{(a-\de)(b-\de)(c-\de)}{\de}\De_\de^++
\tfrac{abc}{\de}+\al(a+b+c+1-d-e)\De_\al^+=0.\label{z42}\end{eqnarray}
In the next 2 formulas $(\de,\ep)$ is a permutation of $(d,e)$.
\begin{eqnarray}
&&\tfrac{(a-\de)(b-\de)(c-\de)}{\de}\De_\de^++
\tfrac{abc}{\de}+(\ep-1)(a+b+c+1-d-e)\De_\ep^-=0.\label{z43}\end{eqnarray}
We conclude by 5 important relations. In first three of them
$(\al,\be,\ga)$ is a cyclic permutation of $(a,b,c)$.
\begin{equation}
(\al-d)(\al-e)\De_\al^-+\al(a+b+c+1-d-e)\De_\al^++
\be\ga=0.\label{z44}\end{equation}
In the last two equalities $\de\in\{d,e\}$.
\begin{equation}
\tfrac{(a-\de)(b-\de)(c-\de)}{\de}\De_\de^++
\tfrac{abc}{\de}+(\de-1)(a+b+c+1-d-e)\De_\de^-=0.\label{z45}\end{equation}

The last 5 relations have the following interpretation:
they are second order difference equations for the function
$F$ considered as a function of the corresponding variable ($a$, $b$, $c$,
$d$ or $e$).

The formulas for the Charlier, Krawtchouk, Meixner, and
Hahn polynomials can be obtained through the corresponding specifications
of the parameters $a,b,c,d,$ and $e$.

Combining some of the above relations (\ref{1})--(\ref{z43})
one can get the following {\em operator shift actions\,}
(the action of the first order difference operators w.r.t. the
variable $c$ is equivalent to the shifting of the parameter $u$):
\begin{eqnarray}
&&\left[\tfrac12(u-\tfrac12)^2+(-c+\tfrac14+\tfrac{d+e-a}{2})(u-\tfrac12)-
\tfrac12(a+\tfrac12)^2+\tfrac12(-de+d+e)\right.\nonumber\\
&&\left.+c(a-\tfrac12)-(c-d)(c-e)\De_c^-+\tfrac{\de}{u-\tfrac12}\right]
\;{}_3F{}_2(a+u,a-u,c;d,e;1)\nonumber\\
&&\qquad =\tfrac{(a-u)(a-d+u)(a-e+u)}{2u-1}\;{}_3F{}_2(a+u-1,a-u+1,c;d,e;1),
\label{n12}\end{eqnarray}
\begin{eqnarray}
&&\left[-\tfrac12(u+\tfrac12)^2+(-c+\tfrac14+\tfrac{d+e-a}{2})(u+\tfrac12)+
\tfrac12(a+\tfrac12)^2-\tfrac12(-de+d+e)\right.\nonumber\\
&&\left.-c(a-\tfrac12)+(c-d)(c-e)\De_c^-+\tfrac{\de}{u+\tfrac12}\right]
\;{}_3F{}_2(a+u,a-u,c;d,e;1)\nonumber\\
&&\qquad =\tfrac{(a+u)(a-d-u)(a-e-u)}{2u+1}\;{}_3F{}_2(a+u+1,a-u-1,c;d,e;1),
\label{n13}\end{eqnarray}
$$
\de=\tfrac12(a-\tfrac12)((a+\tfrac12)(a+\tfrac12-d-e)+de).
$$
The analogous pair is with $c$ shifted in an opposite direction:
\begin{eqnarray}
&&\left[-\tfrac12(u-\tfrac12)^2+(-c-\tfrac34+\tfrac{d+e-a}{2})(u-\tfrac12)+
\tfrac12(a-\tfrac12)^2-\tfrac12(de-d-e+1)\right.\nonumber\\
&&\left.+c(a-\tfrac12)+c(c+2a-d-e+1)\De_c^++\tfrac{\de}{u-\tfrac12}\right]
\;{}_3F{}_2(a+u,a-u,c;d,e;1)\nonumber\\
&&\qquad =\tfrac{(a-u)(a-d+u)(a-e+u)}{2u-1}\;{}_3F{}_2(a+u-1,a-u+1,c;d,e;1),
\label{n14}\end{eqnarray}
\begin{eqnarray}
&&\left[\tfrac12(u+\tfrac12)^2+(-c-\tfrac34+\tfrac{d+e-a}{2})(u+\tfrac12)-
\tfrac12(a-\tfrac12)^2+\tfrac12(de-d-e+1)\right.\nonumber\\
&&\left.-c(a-\tfrac12)-c(c+2a-d-e+1)\De_c^++\tfrac{\de}{u+\tfrac12}\right]
\;{}_3F{}_2(a+u,a-u,c;d,e;1)\nonumber\\
&&\qquad =\tfrac{(a+u)(a-d-u)(a-e-u)}{2u+1}\;{}_3F{}_2(a+u+1,a-u-1,c;d,e;1),
\label{n15}\end{eqnarray}
$$
\de=\tfrac{1}{16}(2a-1)(2a+1-2e)(2a+1-2d).
$$
%
%
\section{Rank 1 quadratic algebra}
\setcounter{equation}{0}
\noindent
The (unitary) $U$-algebra is the algebra with the matrix elements of $U(u)$
as generators and with the quadratic relations given in the form of the
reflection equation (\ref{o4}) (\cite{kk93,ks92,sk88}).
We consider here the reflection equation with the spectral parameter,
with 2-dimensional auxiliary space and with simplest $R$-matrix
of the $XXX$-type (\ref{nn}).

Let us pass to a quotient algebra by adding relations stating that
$(u-\frac12)U(u)$ is a polynomial of degree $\le3$ in $u$.
In other words, we make the ansatz that $U(u)$ is of the form
\begin{eqnarray}
{U(u)}&=&\left(\matrix{ A(u)&B(u)\cr C(u)&D(u)}\right)\nonumber\\
&=&(u-\tfrac12)^2U_2+
(u-\tfrac12)U_1+U_0+
(u-\tfrac12)^{-1}U_{-1},\label{51}\\
\noalign{\hbox{where}}
{U_2}&=&\left(\matrix{ A_2&B_2\cr C_2&-A_2}\right),
\quad U_1=\left(\matrix{ A_1& B_2\cr
 C_2&A_1-2 A_2}\right),\nonumber\\
{U_0}&=&\left(\matrix{ A_0&B_0+\tfrac14B_2\cr
C_0+\tfrac14C_2&-A_0+2 A_1-2A_2}\right),
\quad U_{-1}=\left(\matrix{ A_{-1}&0\cr
0&A_{-1}}\right).\label{51-xx}\end{eqnarray}
Then we get the algebra with the matrix elements of the $U_i$ ($i=2,1,0,-1$)
as generators and with relations
\begin{eqnarray}
{U_2^{(1)}U_{2,1,0,-1}^{(2)}}&=&U_{2,1,0,-1}^{(2)}U_2^{(1)},\quad
U_{-1}^{(1)}U_{2,1,0,-1}^{(2)}=U_{2,1,0,-1}^{(2)}U_{-1}^{(1)},
\label{52}\\
{[U_1^{(1)},U_0^{(2)}]}&=&- [P,U_2^{(1)}U_0^{(2)}]-
 U_2^{(1)}P U_0^{(2)}+ U_0^{(2)}P U_2^{(1)},
\label{53}\\
{[U_0^{(1)},U_0^{(2)}]}&=&- [\{P,U_1^{(1)}\},U_0^{(2)}]-
2 A_{-1}[P,U_2^{(1)}]+{[U_0,U_2]}^{(2)}.
\label{54}\end{eqnarray}
Here curved brackets mean anticommutator.
The relations (\ref{52}) imply that the entries of the $U_2$ and $U_{-1}$
matrices
are in the center of the algebra. Let us pass once more to a quotient algebra
by adding the relations
\begin{equation}
U_2=\left(\matrix{ \alpha&\beta\cr \gamma&-\alpha}\right),
\quad A_{-1}=\delta,
\label{55}\end{equation}
for certain $\al,\be,\ga,\de\in\CC$.
We thus obtain an algebra with generators $A_0,A_1,B_0,C_0$ and relations
\begin{eqnarray}
{[A_1,A_0]}&=&\gamma B_0-\beta C_0,\nonumber\\
{[A_1,B_0]}&=&-2\alpha B_0+\beta\left(2A_0-2A_1+\tfrac32\alpha
\right),\nonumber\\
{[A_1,C_0]}&=&2\alpha C_0+\gamma\left(-2A_0+2A_1-\tfrac32\alpha
\right),\nonumber\\
{[A_0,B_0]}&=&-\{A_1,B_0\}+\beta\left(2A_0-\tfrac52 A_1+2\alpha
+2\delta\right),\nonumber\\
{[A_0,C_0]}&=&\{A_1,C_0\}+\gamma\left(-2A_0+\tfrac52 A_1-2\alpha
-2\delta\right),\nonumber\\
{[B_0,C_0]}&=&-2\{A_0,A_1\}+4(A_1-\alpha)^2+4\alpha A_0+4\alpha\delta\,.
\label{56}\end{eqnarray}
The quantum determinant has the form
\begin{eqnarray}
\Delta(u)&=&-({\alpha}^2+\beta\gamma)u^4+Q_2u^2+Q_0+
{\delta}^2\, u^{-2}\,,\label{57}\\
Q_2&=&A_1^2-2\alpha A_0-\gamma B_0-\beta C_0+\tfrac12\beta\gamma,
\label{t20}\\
Q_0&=&-A_0^2-B_0C_0+2\delta A_1-\tfrac14\gamma B_0-\tfrac14\beta C_0-
\tfrac1{16}\beta\gamma.
\label{t21}\end{eqnarray}
Here the right hand sides of (\ref{t20}) and (\ref{t21}) give operators
in the center of the algebra. We consider (\ref{t20}) and (\ref{t21}) as
added relations, for a certain choice of $Q_2,Q_0\in\CC$.
So a representation of the algebra with relations (\ref{52})--(\ref{54})
which has the
property that all elements in the center of the algebra are represented as
scalars, can also be viewed as a representation of the algebra with
relations (\ref{56}), (\ref{t20}) and (\ref{t21})
for a certain choice of $\al,\be,\ga,\de,Q_0,Q_2$.

Let us change the generators $A_1,A_0,B_0,C_0$ to the new ones
$\taa,\ta,\tb,\tc$:
$$
\taa=A_1-\al,\quad \ta=A_0-A_1+\al,\quad \tb=B_0+\tfrac{\be}{4},
\quad \tc=C_0+\tfrac{\ga}{4}.
$$
The relations for new generators can be rewritten in the following form:
\begin{eqnarray}
\ta\taa&=&\taa\ta-\ga\tb+\be\tc,\nonumber\\
\tb\taa&=&\taa\tb+2\al\tb-2\be\ta,\nonumber\\
\tc\taa&=&\taa\tc-2\al\tc+2\ga\ta,\nonumber\\
\tb\ta&=&\ta\tb+2\taa\tb+2\al\tb-2\be\ta-2\be\de,\nonumber\\
\tc\ta&=&\ta\tc-2\taa\tc+2\al\tc-2\ga\ta+2\ga\de,\nonumber\\
\tc\tb&=&\tb\tc+4\taa\ta-2\ga\tb+2\be\tc-4\al\de.
\label{kkk}\end{eqnarray}
Let we denote by the ${\tilde\FSA}$ the algebra given by the generators
$\taa,\ta,\tb,\tc$ and the relations (\ref{kkk}).
We have the following

\begin{theorem}\label{k7} (Poincare-Birkhof-Witt property)
Let ${\tilde\FSA}^{(n)}\subset{\tilde\FSA}$ be the linear span
of monomials of degree $n$ on generators $\taa,\ta,\tb,\tc$.
Then the dimension of ${\tilde\FSA}^{(n)}$ is equal to the dimension
of the space of monomials of degree $n$ on 4 commuting variables
$\taa,\ta,\tb,\tc$.
\end{theorem}
\vskip 2mm
\Proof
The proof can be done using the Diamond Lemma \cite{be78}.
The defining relations (\ref{kkk}) respect the following ordering
$\taa<\ta<\tb<\tc$. It can be simply verified that there are no
additional relations appearing in the following ambiguities:
$$
\tc\tb\ta,\quad \tc\tb\taa,\quad \tb\ta\taa,\quad \tc\ta\taa.
$$
\vskip 2mm

The group GL(2,$\CC$) acts naturally on the space of $U$-operators
\begin{equation}
{\tilde U}(u)=\left(\matrix{ a&b \cr c&d}\right) U(u)
\left(\matrix{ d&-b \cr -c&a}\right),
\label{n1}\end{equation}
where $ad-bc\neq 0$. We remark that such transformation is compatible to the
unitarity property (\ref{141}).
One may use this transformation to get a more
suitable set of the scalars $\al,\be,\ga$. It is always possible to
arrange that
\begin{equation}
\be=0.
\label{n2}\end{equation}
Then there are only three cases to study:
\begin{eqnarray}
&& Case \quad i):\quad \ga=1,\quad \al=\tfrac12;\nonumber\\
&& Case \quad ii):\quad \ga=1,\quad \al=0;\label{n3}\\
&& Case \quad iii):\quad \ga=\al=0.
\nonumber\end{eqnarray}
The case $ii)$ has been studied in \cite{kk93}.
Here we are dealing more with the cases $i)$ and $iii)$.
%
%
\section{Homomorphisms into $\FSU(g)$}
\setcounter{equation}{0}
\noindent
In this Section we are giving some homomorphisms of the quadratic
algebra with the generators $A_1,A_0,B_0,C_0$ and the relations
(\ref{56}) into the universal enveloping algebra $\FSU(g)$ of the Lie algebra
$g$ for all the cases (\ref{n3}).
%
%
\subsection{The case i)}
Consider the algebra $\FSA$ with generators $A_1,A_0,B_0,C_0,\de$ and with
two sets of relations: relations (\ref{56}) with $\al=0$,
$\be=\ga=1$, and
relations stating that $\de$ is in the center of the algebra.
The case of the scalars $\al,\be,\ga$ just chosen is equivalent
to the case i) in (\ref{n3}) under a certain transformation of the form
(\ref{n1}).
There is a homomorphism of this algebra into the universal enveloping
algebra $\FSU({\rm o}(4))$
of the Lie algebra o(4). A Lie group corresponding to o(4) is the
group of rotations of 4-dimensional Euclidean space.
The Lie algebra o(4) is 6-dimensional.
It can be described by the generators
${\cal P}^\pm,{\cal P}^3,{\cal J}^\pm,{\cal J}^3$ and
commutation relations
$$
[{\cal J}^3,{\cal J}^\pm]=\pm{\cal J}^\pm,\qquad
[{\cal J}^3,{\cal P}^\pm]=[{\cal P}^3,{\cal J}^\pm]=\pm{\cal P}^\pm,
$$
$$
[{\cal J}^+,{\cal P}^+]=[{\cal J}^-,{\cal P}^-]=[{\cal J}^3,{\cal P}^3]=0,
$$
$$
[{\cal J}^+,{\cal J}^-]=2{\cal J}^3,\qquad
[{\cal J}^+,{\cal P}^-]=[{\cal P}^+,{\cal J}^-]=2{\cal P}^3,
$$
$$
[{\cal P}^3,{\cal P}^\pm]=\pm{\cal J}^\pm,\qquad
[{\cal P}^+,{\cal P}^-]=2{\cal J}^3.
$$
The center of the universal enveloping algebra is generated by two
Casimir elements:
$$
C={\left({\cal P}^3\right)}^2+\tfrac12\{{\cal P}^+,{\cal P}^-\}
+{\left({\cal J}^3\right)}^2+\tfrac12\{{\cal J}^+,{\cal J}^-\},
\quad {\tilde C}=\tfrac12\left({\cal P}^+{\cal J}^-+
{\cal P}^-{\cal J}^+\right)+{\cal P}^3{\cal J}^3.
$$
It is now straightforward to verify that the relations for the generators of
$\FSA$ are satisfied when we put these generators equal to the following
elements of $\FSU({\rm o}(4))$.
\begin{eqnarray}
A_1&=&{\cal P}^3,\quad A_0=\tfrac12({\cal P}^+{\cal J}^--{\cal P}^-{\cal J}^+),
\quad B_0=-{\left({\cal J}^3\right)}^2-
\tfrac12\{{\cal P}^+,{\cal P}^-\}-\tfrac14,\nonumber\\
C_0&=&-{\left({\cal J}^3\right)}^2-\tfrac12\{{\cal J}^+,{\cal J}^-\}
-\tfrac14,\quad
\de=-{\tilde C}{\cal J}^3.\qquad\nonumber
\end{eqnarray}
%
%
\subsection{The case ii)}
\noindent
The case ii) ($\al=\be=0,\,\ga=1$) corresponds to the contraction
of the algebra o(4) to the Lie algebra e(3).
So, there is a homomorphism of such algebra into the universal enveloping
algebra $\FSU({\rm e}(3))$.
A Lie group corresponding to e(3) is the
group of motions of 3-dimensional Euclidean space.
The Lie algebra e(3) is 6-dimensional.
It can be described by the generators
${\cal P}^\pm,{\cal P}^3,{\cal J}^\pm,{\cal J}^3$ and
commutation relations
$$
[{\cal J}^3,{\cal J}^\pm]=\pm{\cal J}^\pm,\qquad
[{\cal J}^3,{\cal P}^\pm]=[{\cal P}^3,{\cal J}^\pm]=\pm{\cal P}^\pm,
$$
$$
[{\cal J}^+,{\cal P}^+]=[{\cal J}^-,{\cal P}^-]=[{\cal J}^3,{\cal P}^3]=0,
$$
$$
[{\cal J}^+,{\cal J}^-]=2{\cal J}^3,\qquad
[{\cal J}^+,{\cal P}^-]=[{\cal P}^+,{\cal J}^-]=2{\cal P}^3,
$$
$$
[{\cal P}^3,{\cal P}^\pm]=[{\cal P}^+,{\cal P}^-]=0.
$$
The center of the universal enveloping algebra is generated by two
Casimir elements:
$$
C={\left({\cal P}^3\right)}^2+{\cal P}^+{\cal P}^-,
\qquad {\tilde C}=\tfrac12\left({\cal P}^+{\cal J}^-+
{\cal P}^-{\cal J}^+\right)+{\cal P}^3{\cal J}^3.
$$
The explicit homomorphism has the form \cite{kk93}
\begin{eqnarray}
A_1={\cal P}^3,\quad A_0=\tfrac12({\cal P}^+{\cal J}^--{\cal P}^-{\cal J}^+),
\quad B_0=-{\cal P}^+{\cal P}^-,\nonumber\\
C_0=-{\left({\cal J}^3\right)}^2-\tfrac12\{{\cal J}^+,{\cal J}^-\}
-\tfrac14,\quad
\de=-{\tilde C}{\cal J}^3.\qquad\nonumber
\end{eqnarray}
%
%
\subsection{The case iii)}
\noindent
In the case iii) ($\al=\be=\ga=0$) the generator $A_1$ is in the center
of the algebra and we put $A_1=-\tfrac12$. The commutation relations
(\ref{56}) become linear and we have the following homomorphism of such
algebra to the Lie algebra o(3):
$$
A_0={\cal J}^3-\tfrac12,\quad B_0={\cal J}^+,\quad
C_0={\cal J}^-,\quad
\de\in\CC,
$$
where the o(3) generators
${\cal J}^\pm,{\cal J}^3$ satisfy the commutation relations
$$
[{\cal J}^3,{\cal J}^\pm]=\pm{\cal J}^\pm,\quad
[{\cal J}^+,{\cal J}^-]=2{\cal J}^3.
$$
%
%
\section{Realisation via difference operators}
\setcounter{equation}{0}
\noindent
Let us consider the generic case i) ($\be=0$, $\al=\tfrac12,\ga=1$).
The following lemma is a slightly extended form of the Lemma~6.2
in \cite{kk93}. It
shows that equations (\ref{56}), (\ref{t20}) and (\ref{t21}), with
(\ref{n2}) ( and with (\ref{n3}) for the case i) ),
and under the assumption that $B_0$ is injective,
can be equivalently written in a much more simple form.

\begin{lemma}\label{t23}
Let $\delta,Q_0,Q_2$ be scalars.
Let $A_1,A_0,B_0,C_0$
be operators acting on some linear space. Let $B_0$ be injective.
Then the following three statements are equivalent:
\vskip 2mm\par\noindent
(a)
$\;\left(\matrix{ \tfrac12(u-\tfrac12)^2+(u-\tfrac12)
A_1+A_0+\tfrac{\de}{u-\tfrac12}&B_0\cr
{}&-\tfrac12(u+\tfrac12)^2+(u+\tfrac32)
A_1\cr u^2+C_0&
-A_0-\tfrac12+\tfrac{\de}{u-\tfrac12}}\right)$
\vskip 2mm\par\noindent
is a representation of the $U$-algebra with quantum determinant
$\De(u)=-\tfrac14u^4+Q_2u^2+Q_0+\de^2u^{-2}$;
\vskip 1mm

\noindent
(b)\quad
The six commutators (\ref{56}) and formulas (\ref{t20}), (\ref{t21})
are valid with $\al=\tfrac12,\be=0$ and $\ga=1$;
\vskip 1mm

\noindent
(c)\quad
The following three equations are valid:
\begin{eqnarray}
[A_1,A_0]&=&A_1^2-A_0-Q_2,\label{59}\\
{B_0}&=&A_1^2-A_0-Q_2,\label{58}\\
B_0C_0&=&2\de A_1-A_0^2-\tfrac14  B_0-Q_0.\label{t24}
\end{eqnarray}
\end{lemma}

We now make the restrictive assumption
that $A_0$ is a second order difference
operator and $A_1$ is a scalar function of $x$:
\begin{equation}
A_0=A_0^-(x)\;(1+\De_x^-)+A_0^+(x)\;(1+\De_x^+)+A_{00}(x),\quad A_1=A_{10}(x).
\label{510}
\end{equation}
The following approach should now be followed.
Find all operators $A_0,A_1$ of the form (\ref{510})
such that (\ref{59}) is satisfied for some number $Q_2$.
(It is sufficient to find one solution in each equivalence class formed
by gauge transformations.)$\;$
Then define $B_0$ by (\ref{58}) and try to define $C_0$ by (\ref{t24}).
Fix some function space $W$ on which these operators act.
Then the equivalent conditions of Lemma \ref{t23} are satisfied.
Finally check if the conditions of Theorem \ref{t7} are satisfied for some
choice of $\FSD$.

Below we are giving two lemmas which can be proved by
straightforward computation.

\begin{lemma}\label{t60}
If we assume (\ref{510}) for the $A_0$ and $A_1$ then (\ref{59})
holds if and only if
\begin{eqnarray}
&&A_{00}(x)={A_{10}}^2(x)-Q_2,\label{kkk1}\\
&&A_0^-(x)\,\left(A_{10}(x)-A_{10}(x-1)+1\right)=0,\label{kkk2}\\
&&A_0^+(x)\,\left(A_{10}(x)-A_{10}(x+1)+1\right)=0.\label{kkk3}\end{eqnarray}
\end{lemma}

\begin{lemma}\label{tt60}
Let two functions $A_0^\pm(x)$ be not both identically equal to zero.
Then there are only two solutions of the equations (\ref{kkk1})--(\ref{kkk3}):
\begin{eqnarray}
&&(a)
\qquad A_0^+(x)\equiv 0,\quad A_{10}(x)=-x+c_1,\quad A_{00}(x)=(x-c_1)^2-Q_2,
\nonumber\\
&&(b)
\qquad A_0^-(x)\equiv 0,\quad A_{10}(x)=x+c_1,\quad A_{00}(x)=(x+c_1)^2-Q_2
\nonumber\end{eqnarray}
where $c_1$ is an arbitrary constant and one has an arbitrary
function ( $A_0^-(x)$ in the case (a) and $A_0^+(x)$ in the case (b) )
which can be fixed by the gauge transformation.
\end{lemma}
\vskip 2mm

The case $(a)$ of the Lemma \ref{tt60} is equivalent (up to a gauge
transformation and a constant shift of the independent variable $x$)
to the following realisation of the algebra in terms of the difference
operators:
\begin{eqnarray}
A_1&=&-x+\tfrac14+\tfrac{d+e-a}{2}\,,\nonumber\\
A_0&=&-\tfrac12(a+\tfrac12)^2+\tfrac12(-de+d+e)
+x(a-\tfrac12)-(x-d)(x-e)\De_x^-,\nonumber\\
B_0&=&(x-d)(x-e)(1+\De_x^-),\nonumber\\
C_0&=&-(x-d)(x-e)\De_x^--x(x+2a+1-d-e)\De_x^+-a^2,\nonumber\\
Q_2&=&{\tfrac {a}{4}}-{\tfrac {d}{4}}-{\tfrac {e}{4}}+{\tfrac {3\,a^{2}}{4}}-{
\tfrac {ad}{2}}-{\tfrac {ae}{2}}+{\tfrac {d^{2}}{4}}+{\tfrac {e^{2}}{4}}+3
/16\,,\nonumber\\
Q_0&=&-{\tfrac {ed}{4}}-{\tfrac {a}{8}}+{\tfrac {d}{8}}+{\tfrac {e}{8}}-eda-{
\tfrac {a^{2}}{8}}+{\tfrac {ad}{4}}+{\tfrac {ae}{4}}-{\tfrac {d^{2}}{8}}
\nonumber\\
&&-{\tfrac {e^{2}}{8}}+{\tfrac {ed^{2}}{4}}+{\tfrac {e^{2}d}{4}}
+{\tfrac {ed^{2}a}{2}}+{\tfrac {e^{2}da}{2}}-{\tfrac {3}{64}}-
{\tfrac {3\,a^{4}}{4}}\nonumber\\
&&+{\tfrac {a^{2}d}{2}}+{\tfrac {a^{2}e}{2}}-{\tfrac {a^{3}}{2}}+a^{3}d+a^{3}
e-{\tfrac {e^{2}d^{2}}{4}}-{\tfrac {a^{2}d^{2}}{2}}-{\tfrac {a^{2}e^{2}}{
2}}-a^{2}ed\,,\nonumber\\
\de&=&\tfrac12(a-\tfrac12)((a+\tfrac12)(a+\tfrac12-d-e)+de).
\nonumber\end{eqnarray}

The case $(b)$ of the Lemma \ref{tt60} is equivalent (up to a gauge
transformation and a constant shift of the independent variable $x$)
to the following realisation of the algebra in terms of the difference
operators:
\begin{eqnarray}
A_1&=&x+\tfrac34-\tfrac{d+e-a}{2}\,,\nonumber\\
A_0&=&-\tfrac12(a-\tfrac12)^2+\tfrac12(de-d-e+1)-x(a-\tfrac12)-
x(x+2a+1-d-e)\De_x^+,\nonumber\\
B_0&=&x(x+2a+1-d-e)(1+\De_x^+),\nonumber\\
C_0&=&-(x-d)(x-e)\De_x^--x(x+2a+1-d-e)\De_x^+-a^2,\nonumber\\
Q_2&=&{\tfrac {a}{4}}-{\tfrac {e}{4}}-{\tfrac {d}{4}}+{\tfrac {3\,a^{2}}{4}}+{
\tfrac {e^{2}}{4}}-{\tfrac {ea}{2}}+{\tfrac {d^{2}}{4}}+3/16-{\tfrac {da}{
2}}\,,\nonumber\\
Q_0&=&
-{\tfrac {de}{4}}-{\tfrac {a}{8}}+{\tfrac {e}{8}}+{\tfrac {d}{8}}-{\tfrac {
a^{2}}{8}}-dea-{\tfrac {e^{2}}{8}}+{\tfrac {ea}{4}}+{\tfrac {da}{4}}-{
\tfrac {d^{2}}{8}}\nonumber\\
&&+{\tfrac {de^{2}}{4}}+{\tfrac {d^{2}e}{4}}-a^{2}de+{\tfrac {de^{2}a}{2}}
+{\tfrac {d^{2}ea}{2}}-{\tfrac {3\,a^{4}}{4}}-{\tfrac {3}{64}}+a^{3}e+a^{3}d
\nonumber\\
&&-{\tfrac {a^{2}d^{2}}{2}}-{\tfrac {a^{2}e^{2}}{2}}-{\tfrac {d^{2}e^{2}}{4}}
+{\tfrac {a^{2}e}{2}}+{\tfrac {a^{2}d}{2}}-{\tfrac {a^{3}}{2}}\,,\nonumber\\
\de&=&\tfrac{1}{16}(2a-1)(2a+1-2e)(2a+1-2d).
\nonumber\end{eqnarray}
Here the $a,d,$ and $e$ are arbitrary parameters.

These operators act on the functions ${}_3F{}_2(a+u,a-u,x;d,e;1)$
(cf. (3.13)--(3.14) and (3.15)--(3.16)).
%
%
\section{Concluding remarks}
\setcounter{equation}{0}
\noindent
It follows from (\ref{o10}) (second and fourth formula) that
\begin{eqnarray}
-A(-u+1)A(u)&=&B(u-1)C(u)+\De(u-\tfrac12),\label{xt34}\\
-A(u+1)A(-u)&=&B(u+1)C(u)+\De(u+\tfrac12).\label{xt35}
\end{eqnarray}
So assume that we have a representation of the QISM II algebra on a space
of functions in one variable, analytic on a certain region, such that
(i) $B(u)=B_0$ is independent of $u$,
(ii) the unitarity property (\ref{141}) is satisfied,
(iii) $\De(u)$ is scalar for all $u$,
(iv) $A(u)$ is a first order difference operator.
Then equations (\ref{xt34})--(\ref{xt35}) show that the second order
difference operator
$B_0C(u)$ has a suitable form for the factorization method,
which was originated by Schr\"odinger and was due in its definitive form to
Infeld and Hull \cite{infeld}.

The factorization method for the second order differential operators
under some special assumptions on the type of factorizing operators
has been summarized in Miller \cite[Ch.\ 7]{mi68}.
These factorizing operators give rise to a Lie algebra of
first order differential operators in two variables or to some operators
in the universal enveloping algebra of certain Lie algebras.
See \cite{kk93} as for the $R$-matrix interpretation.
The generalization to the difference operators
and some further generalizations have been given in \cite{mi69,mi70}.

In this work we have established a connection of the QISM II algebra
representations to some characterizing properties of the generalized
hypergeometric function ${}_3F{}_2(1)$. Namely, we have given
the algebraic interpretation of the three relations between
${}_3F{}_2(1)$ and some other two functions which are contiguous to it. On the
Hahn polynomials level these formulas are equivalent to the
recurrence relation, difference ``differentiation'' formula and
the second order difference equation for the polynomials
(operators $A_1, A_0$, and $C_0$).

The presented homomorphisms of the QISM II algebra into $\FSU(g)$
will allow us to construct some new integrable systems which will be
published elsewhere.

\bibliographystyle{amsplain}

\begin{thebibliography}{10}

\bibitem{ba82}
R.J.~Baxter.
\newblock {\em Exactly solved models in statistical mechanics}.
\newblock Academic Press, New-York, 1982.

\bibitem{be78}
G.M.~Bergman.
\newblock The diamond lemma for ring theory.
\newblock {\em Advances Math.} 29:178--218, 1978.

\bibitem{fa94}
L.D.~Faddeev.
\newblock Algebraic aspects of Bethe-ansatz.
\newblock Preprint, Stony Brook, ITP-SB-94-11, 42pp., March 1994.
          HEP-TH/9404013

\bibitem{infeld}
L.~Infeld and T.~Hull.
\newblock The factorization method.
\newblock {\em Rev. Mod. Phys.} 23:21--68, 1951.

\bibitem{kk93}
T.H.~Koornwinder and V.B.~Kuznetsov.
\newblock Gauss hypergeometric function and quadratic ${R}$-matrix algebras.
\newblock {\em Algebra and Analysis} 6:161--184, 1994.

\bibitem{ks92}
P.P.~Kulish and E.K.~Sklyanin.
\newblock Algebraic structures related to reflection equations.
\newblock {\em J. Phys. A: Math.Gen.} 25:5963--5975, 1992.

\bibitem{mi68}
W.~Miller{,} Jr.
\newblock {\em Lie theory and special functions}.
\newblock Academic Press, New{-}York, 1968.

\bibitem{mi69}
W.~Miller{,} Jr.
\newblock Lie theory and difference equations, I.
\newblock {\em J. Math. Anal. Appl.} 28:383--399, 1969.

\bibitem{mi70}
W.~Miller{,} Jr.
\newblock Lie theory and q-difference equations.
\newblock {\em SIAM J. Math. Anal.} 1:171-188, 1970.

\bibitem{ra45}
E.D.~Rainville.
\newblock The contiguous function relations for ${}_pF_q$ with applications
          to Bateman's $J_n^{u,v}$ and Rice's $H_n(z,p,v)$.
\newblock {\em Bull. Amer. Math. Soc.}, 51:714--723, 1945.

\bibitem{sk88}
E.K.~Sklyanin.
\newblock Boundary conditions for quantum integrable systems.
\newblock {\em J. Phys. A}, 21:2375--2387, 1988.

\bibitem{sk91}
E.K.~Sklyanin.
\newblock Quantum inverse scattering method. {S}elected topics.
\newblock In: Mo-Lin Ge, editor, {\em {Q}uantum {G}roup and {Q}uantum
{I}ntegrable {S}ystems}, (Nankai Lectures in Mathematical Physics),
Singapore: World Scientific, 63--97, 1992.

\end{thebibliography}

\end{document}